\documentclass{article}
\usepackage{amsfonts}
\usepackage{amssymb}
\usepackage{graphicx}
\usepackage{amsmath}

\input{tcilatex}
\begin{document}

\begin{center}
\textbf{\ Mixed problems for degenerate abstract parabolic equations and
applications}

\ 

\textbf{Veli.B. Shakhmurov}

Department of Mechanical Engineering, Okan University, Akfirat, Tuzla 34959
Istanbul, Turkey,

E-mail: veli.sahmurov@okan.edu.tr\ 

Khazar University, Baku Azerbaijan

\textbf{Aida}{\textbf{\ Sahmurova}}

Okan University, Department of Environmental Engineering, Akfirat, Tuzla
34959 Istanbul, Turkey, E-mail: aida.sahmurova@okan.edu.tr
\end{center}


\begin{center}
\textbf{AMC 2000: 35A01, 35J56, 35Dxx, 35K51, 47G40}

\bigskip\ \ \ \ \ \ \ \ \ \ \ 

\textbf{ABSTRACT}
\end{center}

Degenerate abstract parabolic equations with variable coefficients are
studied. Here the boundary conditions are nonlocal. The maximal regularity
properties of solutions for elliptic and parabolic problems and Strichartz
type estimates in mixed $L_{\mathbf{p}}$ spaces\ are obtained. Moreover, the
existence and uniqueness of optimal regular solution of mixed problem for
nonlinear parabolic equation is established. Note that, these problems arise
in fluid mechanics and environmental engineering.

\begin{center}
\ \ \textbf{1. Introduction and notations }
\end{center}

In this work, the boundary value problems (BVPs) for parameter dependent
degenerate differential-operator equations (DOEs)\ are considered. Namely,
linear equations and boundary conditions contain small parameters and are
degenerated in some part of boundary. These problems have numerous
applications in PDE, pseudo DE, mechanics and environmental engineering. The
BVP for DOEs have been studied extensively by many researchers (see e.g. $%
\left[ \text{1-11}\right] $ and the references therein). The maximal
regularity properties for DOEs in Banach space valued function class are
investigated e.g. in $\left[ \text{4-11}\right] .$ Nonlinear DOEs studied
e.g. in $\left[ \text{7,10}\right] $.

The main objective of the present paper is to discuses the initial and
nonlocal BVP for the following nonlinear degenerate parabolic equation 
\begin{equation}
\frac{\partial u}{\partial t}+\dsum\limits_{k=1}^{n}a_{k}\left( x\right) 
\frac{\partial ^{\left[ 2\right] }u}{\partial x_{k}^{2}}+B\left( \left(
t,x,u,D^{\left[ 1\right] }u\right) \right) u=F\left( t,x,u,D^{\left[ 1\right]
}u\right) ,  \tag{1.1}
\end{equation}%
where $a_{k}$ are complex valued functions, $B$ and $F$ are nonlinear\
operators in a Banach space $E$ and 
\begin{equation*}
D^{\left[ 1\right] }u=\left( \frac{\partial ^{\left[ 1\right] }u}{\partial
x_{1}},\frac{\partial ^{\left[ 1\right] }u}{\partial x_{2}},...,\frac{%
\partial ^{\left[ 1\right] }u}{\partial x_{n}}\right) ,\text{ }x=\left(
x_{1},x_{2},...,x_{n}\right) \in G=\dprod\limits_{k=1}^{n}\left(
0,b_{k}\right) \text{,}
\end{equation*}

\begin{equation*}
D_{_{k}}^{\left[ i\right] }u=u_{k}^{\left( i\right) }=\frac{\partial ^{\left[
i\right] }u}{\partial x_{k}^{i}}=\left( x_{k}^{\alpha _{k}}\frac{\partial }{%
\partial x_{k}}\right) ^{i}u\left( x\right) ,\text{ }0\leq \alpha _{k}<1.
\end{equation*}

First all of, we consider the nonlocal BVP for the degenerate elliptic DOE
with small parameters

\begin{equation}
\dsum\limits_{k=1}^{n}a_{k}\left( x\right) \frac{\partial ^{\left[ 2\right]
}u}{\partial x_{k}^{2}}+A\left( x\right) u+\lambda
u+\dsum\limits_{k=1}^{n}A_{k}\left( x\right) \frac{\partial ^{\left[ 1\right]
}u}{\partial x_{k}}=f\left( x\right) ,  \tag{1.2}
\end{equation}%
where $a_{k}$ are complex-valued functions, $\lambda $ is a complex
parameters, $A\left( x\right) $ and $A_{k}\left( x\right) $ are linear
operators.

We prove that for $f\in L_{p}\left( G;E\right) $, $\left\vert \arg \lambda
\right\vert \leq \varphi ,$ $0<\varphi \leq \pi $ and sufficiently large $%
\left\vert \lambda \right\vert ,$ problem $\left( 1.2\right) $ has a unique
solution $u\in $ $W_{p\mathbf{,\alpha }}^{\left[ 2\right] }\left( G;E\left(
A\right) ,E\right) $ and the following coercive uniform estimate holds

\begin{equation*}
\sum\limits_{k=1}^{n}\sum\limits_{i=0}^{2}\left\vert \lambda \right\vert ^{1-%
\frac{i}{2}}\left\Vert \frac{\partial ^{\left[ i\right] }u}{\partial
x_{k}^{i}}\right\Vert _{L_{p}\left( G;E\right) }+\left\Vert Au\right\Vert
_{L_{p}\left( G;E\right) }\leq C\left\Vert f\right\Vert _{L_{p}\left(
G;E\right) }.
\end{equation*}%
Then the above result is used to prove the\ well-posedeness of initial BVP
(IBVP) and the uniform Strichartz type estimate for the solution the
degenerate abstract parabolic equation with parameters%
\begin{equation}
\frac{\partial u}{\partial t}+\dsum\limits_{k=1}^{n}a_{k}\left( x\right) 
\frac{\partial ^{\left[ 2\right] }u}{\partial x_{k}^{2}}+A\left( x\right)
u=f\left( x,t\right) \text{, }t\in \left( 0,T\right) \text{, }x\in G. 
\tag{1.3}
\end{equation}%
Finally, via maximal regularity properties of $\left( 1.3\right) $ and
contraction mapping argument, the existence and uniqueness of solution of
the problem $\left( 1.1\right) $ is derived.

Note that, the equation and boundary conditions are degenerated with the
different rate at different boundary edges, in general.

In application, the system of degenerate nonlinear parabolic equations is
presented. Particularly, we consider the system that serves as a model of
systems used to describe photochemical generation and atmospheric dispersion
of ozone and other pollutants. The model of the process is given by initial
and BVP for the atmospheric reaction-advection-diffusion system having the
form 
\begin{equation}
\frac{\partial u_{i}}{\partial t}=\sum\limits_{k=1}^{3}\left[ a_{ki}\left(
x\right) \frac{\partial ^{\left[ 2\right] }u_{i}}{\partial x_{k}^{2}}%
+b_{ki}\left( x\right) \frac{\partial ^{\left[ 1\right] }}{\partial x_{k}}%
\left( u_{i}\omega _{k}\right) \right] +\sum%
\limits_{k=1}^{3}d_{k}u_{k}+f_{i}\left( u\right) +g_{i}\text{,}  \tag{0.4}
\end{equation}%
where 
\begin{equation*}
\text{ }x\in G_{3}=\left\{ x=\left( x_{1},x_{2},x_{3}\right) ,\text{ }%
0<x_{k}<b_{k},\right\} \text{, }
\end{equation*}%
\begin{equation*}
u_{i}=u_{i}\left( x,t\right) ,\text{ }i,\text{ }k=1,2,3,\text{ }u=u\left(
x,t\right) =\left( u_{1},u_{2},u_{3}\right) ,\text{ }t\in \left( 0,T\right)
\end{equation*}%
and the state variables $u_{i}$ represent concentration densities of the
chemical species involved in the photochemical reaction. The relevant
chemistry of the chemical species involved in the photochemical reaction and
appears in the nonlinear functions $f_{i}\left( u\right) ,$ with the terms $%
g_{i},$ representing elevated point sources, $a_{ki}\left( x\right)
,b_{ki}\left( x\right) $ are real-valued functions.\ The advection terms $%
\omega =\omega \left( x\right) =\left( \omega _{1}\left( x\right) ,\omega
_{2}\left( x\right) ,\omega _{3}\left( x\right) \right) $, describe
transport from the velocity vector field of atmospheric currents or wind. In
this direction the work $\left[ 12\right] $ and references there can be
mentioned$.$ The existence and uniqueness of solution of the problem $\left(
0.4\right) $ is established by the theoretic-operator method, i.e., this
problem is reduced to degenerate differential-operator equation.

Let $\gamma =\gamma \left( x\right) $ be a positive measurable function on $%
\Omega \subset R^{n}$ and $E$ be a Banach space.\ Let $L_{p,\gamma }\left(
\Omega ;E\right) $ denote the space of strongly measurable $E$-valued
functions defined on $\Omega $ with the norm

\begin{equation*}
\left\Vert f\right\Vert _{L_{p,\gamma }}=\left\Vert f\right\Vert
_{L_{p,\gamma }\left( \Omega ;E\right) }=\left( \int \left\Vert f\left(
x\right) \right\Vert _{E}^{p}\gamma \left( x\right) dx\right) ^{\frac{1}{p}%
},1\leq p<\infty .
\end{equation*}

For $\gamma \left( x\right) \equiv 1$ we will denote these spaces by $%
L_{p}\left( \Omega ;E\right) .$

The Banach space\ $E$ is called an $UMD$-space if\ the Hilbert operator%
\begin{equation*}
\left( Hf\right) \left( x\right) =\lim\limits_{\varepsilon \rightarrow
0}\int\limits_{\left\vert x-y\right\vert >\varepsilon }\frac{f\left(
y\right) }{x-y}dy
\end{equation*}%
\ is bounded in $L_{p}\left( R,E\right) ,$ $p\in \left( 1,\infty \right) $ (
see. e.g. $\left[ 13\right] $ ). $UMD$ spaces include e.g. $L_{p}$, $l_{p}$
spaces and Lorentz spaces $L_{pq},$ $p$, $q\in \left( 1,\infty \right) $.

Let $\mathbb{C}$ be the set of the complex numbers and\ 
\begin{equation*}
S_{\varphi }=\left\{ \lambda ;\text{ \ }\lambda \in \mathbb{C}\text{, }%
\left\vert \arg \lambda \right\vert \leq \varphi \right\} \cup \left\{
0\right\} ,\text{ }0\leq \varphi <\pi .
\end{equation*}

Let $E_{1}$ and $E_{2}$ be two Banach spaces. $L\left( E_{1},E_{2}\right) $
denotes the space of bounded linear operators from $E_{1}$ into $E_{2}.$ For 
$E_{1}=E_{2}=E$ the space $L\left( E_{1},E_{2}\right) $ will be denoted by $%
L\left( E\right) .$ A linear operator\ $A$ is said to be $\varphi $-positive
in a Banach\ space $E$ with bound $M>0$ if $D\left( A\right) $ is dense on $%
E $ and $\left\Vert \left( A+\lambda I\right) ^{-1}\right\Vert _{L\left(
E\right) }\leq M\left( 1+\left\vert \lambda \right\vert \right) ^{-1}$ for
any $\lambda \in S_{\varphi },$ $0\leq \varphi <\pi ,$ where $I$ is the
identity operator in $E.$

Let $E_{0}$ and $E$ be two Banach spaces and $E_{0}$ is continuously and
densely embeds into $E$. Let us consider the Sobolev-Lions type space\ $%
W_{p,\gamma }^{m}\left( a,b;E_{0},E\right) ,$ consisting of all functions $%
u\in L_{p,\gamma }\left( a,b;E_{0}\right) $ that have generalized
derivatives $u^{\left( m\right) }\in L_{p,\gamma }\left( a,b;E\right) $ with
the norm 
\begin{equation*}
\ \left\Vert u\right\Vert _{W_{p,\gamma }^{m}}=\left\Vert u\right\Vert
_{W_{p,\gamma }^{m}\left( a,b;E_{0},E\right) }=\left\Vert u\right\Vert
_{L_{p,\gamma }\left( a,b;E_{0}\right) }+\left\Vert u^{\left( m\right)
}\right\Vert _{L_{p,\gamma }\left( a,b;E\right) }<\infty .
\end{equation*}

Let $\gamma =\gamma \left( x\right) $ be a positive measurable function on $%
\left( 0,1\right) $ and 
\begin{equation*}
W_{p,\gamma }^{\left[ m\right] }=W_{p,\gamma }^{\left[ m\right] }\left(
0,1;E_{0},E\right) =\left\{ u:u\in L_{p}\left( 0,1;E_{0}\right) \right. ,\ 
\end{equation*}

\begin{equation*}
\ u^{\left[ m\right] }\in L_{p}\left( 0,1;E\right) ,\left\Vert u\right\Vert
_{W_{p,\gamma }^{\left[ m\right] }}=\left. \left\Vert u\right\Vert
_{L_{p}\left( 0,1;E_{0}\right) }+\left\Vert u^{\left[ m\right] }\right\Vert
_{L_{p}\left( 0,1;E\right) }<\infty \right\} .
\end{equation*}
Let%
\begin{equation*}
\alpha _{k}\left( x\right) =x_{k}^{\alpha _{k}},\text{ }\alpha =\left(
\alpha _{1},\alpha _{2},...,\alpha _{n}\right) .
\end{equation*}

Consider $E$-valued weighted space defined by 
\begin{equation*}
W_{p\mathbf{,}\alpha }^{\left[ m\right] }\left( G,E\left( A\right) ,E\right)
=\left\{ u;u\in L_{p}\left( G;E_{0}\right) \right. ,\ \frac{\partial ^{\left[
m\right] }u}{\partial x_{k}^{m}}\in L_{p}\left( G;E\right) ,
\end{equation*}

\begin{equation*}
\ \left\Vert u\right\Vert _{W_{p\mathbf{,}\alpha }^{\left[ m\right]
}}=\left. \left\Vert u\right\Vert _{L_{p}\left( G;E_{0}\right)
}+\dsum\limits_{k=1}^{n}\left\Vert \frac{\partial ^{\left[ m\right] }u}{%
\partial x_{k}^{m}}\right\Vert _{L_{p}\left( G;E\right) }<\infty \right\} .
\end{equation*}

\begin{center}
\textbf{2. Degenerate abstract elliptic equations }
\end{center}

Consider the BVP for the following degenerate partial DOE with parameters

\begin{equation}
\dsum\limits_{k=1}^{n}a_{k}\left( x_{k}\right) \frac{\partial ^{\left[ 2%
\right] }u}{\partial x_{k}^{2}}+A\left( x\right) u+\lambda
u+\dsum\limits_{k=1}^{n}A_{k}\left( x\right) \frac{\partial ^{\left[ 1\right]
}u}{\partial x_{k}}=f\left( x\right) ,  \tag{2.1}
\end{equation}

\begin{equation*}
L_{kj}u=\sum\limits_{i=0}^{m_{kj}}\alpha _{kji}u_{x_{k}}^{\left[ i\right]
}\left( G_{k0}\right) +\beta _{kji}u_{k}^{\left[ i\right] }\left(
G_{kb}\right) =0,\text{ }j=1,2,
\end{equation*}%
where $a_{k}$\ are complex-valued functions, $A\left( x\right) $ and $%
A_{k}\left( x\right) $ are linear operators, $u=u\left( x\right) $, $\alpha
_{kji},$ $\beta _{kji}$ are complex numbers, $\lambda $ is a complex
parameter, $m_{kj}\in \left\{ 0,1\right\} $, 
\begin{equation*}
x=\left( x_{1},x_{2},...,x_{n}\right) \in G=\dprod\limits_{k=1}^{n}\left(
0,b_{k}\right) \text{, }
\end{equation*}%
\begin{equation*}
\text{ }G_{k0}=\left( x_{1},x_{2},...,x_{k-1},0,x_{k+1},...,x_{n}\right) ,%
\text{ }p_{k}\in \left( 1,\infty \right) ,
\end{equation*}

\begin{equation*}
G_{kb}=\left( x_{1},x_{2},...,x_{k-1},b_{k},x_{k+1},...,x_{n}\right) \text{,}
\end{equation*}%
\begin{equation*}
\text{ }x^{\left( k\right) }=\left(
x_{1},x_{2},...,x_{k-1},x_{k+1},...,x_{n}\right) \in
G_{k}=\dprod\limits_{j\neq k}\left( 0,b_{j}\right) .\text{ }
\end{equation*}

Consider the principal part of $\left( 2.1\right) $, i.e., consider the
problem 
\begin{equation}
\dsum\limits_{k=1}^{n}a_{k}\left( x_{k}\right) \frac{\partial ^{\left[ 2%
\right] }u}{\partial x_{k}^{2}}+A\left( x\right) u+\lambda u=f\left(
x\right) ,  \tag{2.2}
\end{equation}

\begin{equation*}
\sum\limits_{i=0}^{m_{kj}}\alpha _{kji}u_{x_{k}}^{\left[ i\right] }\left(
G_{k0}\right) +\beta _{kji}u_{k}^{\left[ i\right] }\left( G_{kb}\right) =0,%
\text{ }j=1,2.
\end{equation*}

\textbf{Condition 2.1 }Assume;

(1) $E$ is an UMD spacethe Banach space, $0\leq \alpha _{k}<1-\frac{1}{p_{k}}
$, $p_{k}\in \left( 1,\infty \right) $, $\alpha _{km_{k1}}\neq 0$, $\beta
_{km_{k2}}\neq 0;$

(2)\ $A\left( x\right) $ is a uniformly $R$-positive operator in $E$, $%
A\left( x\right) A^{-1}\left( \bar{x}\right) \in C\left( \bar{G};L\left(
E\right) \right) $, $x\in G;$

(3) $a_{k}\in C^{\left( m\right) }\left( \left[ 0,b_{k}\right] \right) $ and 
$a_{k}\left( x_{k}\right) <0$ for $x_{k}\in \left[ 0,b_{k}\right] ;$

(4) $a_{k}\left( G_{j0}\right) =a_{k}\left( G_{jb}\right) ,$ $A\left(
G_{j0}\right) A^{-1}\left( x_{0}\right) =A\left( G_{jb}\right) A^{-1}\left(
x_{0}\right) $, $k,j=1,2,...,n;$

(5) $\eta _{k}=\left( -1\right) ^{m_{1}}\alpha _{k1}\beta _{k2}-\left(
-1\right) ^{m_{2}}\alpha _{k2}\beta _{k1}\neq 0.$

First, we prove the separability properties of the problem $\left(
2.2\right) $:

\textbf{Theorem 2.1. }Let the Conditions 2.1 hold. Then, problem $\left(
2.2\right) $ has a unique solution $u\in W_{p\mathbf{,}\alpha }^{\left[ 2%
\right] }\left( G;E\left( A\right) ,E\right) $ for $f\in L_{p}\left(
G;E\right) $, $\left\vert \arg \lambda \right\vert \leq \varphi $ with
sufficiently large $\left\vert \lambda \right\vert $ and the following
coercive uniform estimate holds

\begin{equation}
\sum\limits_{k=1}^{n}\sum\limits_{i=0}^{2}\left\vert \lambda \right\vert ^{1-%
\frac{i}{2}}\left\Vert \frac{\partial ^{\left[ i\right] }u}{\partial
x_{k}^{i}}\right\Vert _{L_{p}\left( G;E\right) }+\left\Vert Au\right\Vert
_{L_{p}\left( G;E\right) }\leq C\left\Vert f\right\Vert _{L_{p}\left(
G;E\right) }.  \tag{2.3}
\end{equation}

\textbf{Proof. }Consider the BVP 
\begin{equation}
\ \left( L+\lambda \right) u=a_{1}\left( x_{1}\right) D_{x_{1}}^{\left[ 2%
\right] }u\left( x_{1}\right) +\left( A\left( x_{1}\right) +\lambda \right)
u\left( x_{1}\right) =f\left( x_{1}\right) ,\text{ }  \tag{2.4}
\end{equation}

\begin{equation*}
\text{ }L_{1j}u=0,\text{ }j=1,2,\text{ }x_{1}\in \left( 0,b_{1}\right) ,
\end{equation*}%
where $L_{1j}$ are boundary conditions of type $\left( 2.2\right) $
considered on $\left( 0,b_{1}\right) .$ By virtue of Theorem 1 in $\left[ 
\text{8}\right] $, problem $\left( 2.4\right) $ has a unique solution $u\in
W_{p,\alpha _{1}}^{\left[ 2\right] }\left( 0,b_{1};E\left( A\right)
,E\right) $ for $f\in L_{p_{1}}\left( 0,b_{1};E\right) $, $\left\vert \arg
\lambda \right\vert \leq \varphi $ with sufficiently large $\left\vert
\lambda \right\vert $ and the coercive uniform estimate holds

\begin{equation*}
\sum\limits_{j=0}^{2}\left\vert \lambda \right\vert ^{1-\frac{j}{2}%
}\left\Vert u^{\left[ j\right] }\right\Vert _{L_{p_{1}}\left(
0,b_{1};E\right) }+\left\Vert Au\right\Vert _{L_{p_{1}}\left(
0,b_{1};E\right) }\leq C\left\Vert f\right\Vert _{L_{p_{1}}\left(
0,b_{1};E\right) }.
\end{equation*}%
Now, let us consider the following BVP 
\begin{equation}
\sum\limits_{k=1}^{2}a_{k}\left( x_{k}\right) D_{_{k}}^{\left[ 2\right]
}u\left( x_{1},x_{2}\right) +A\left( x_{1},x_{2}\right) u\left(
x_{1},x_{2}\right) +\lambda u\left( x_{1},x_{2}\right) =f\left(
x_{1},x_{2}\right) ,  \tag{2.5}
\end{equation}

\begin{equation*}
L_{k1}u=0,\text{ }L_{k2}u=0,\text{ }k=1,2\text{, }x_{1},x_{2}\in
G_{2}=\left( 0,b_{1}\right) \times \left( 0,b_{2}\right) .
\end{equation*}

Let $\alpha \left( 2\right) =\left( \alpha _{1},\alpha _{2}\right) $. Since $%
L_{p}\left( 0,b_{2};L_{p}\left( 0,b_{1}\right) \text{; }E\right) =$ $%
L_{p}\left( G_{2};E\right) ,$ the BVP $\left( 2.5\right) $ can be expressed
as

\begin{equation*}
a_{2}D_{2}^{\left[ 2\right] }u\left( x_{2}\right) +\left( B\left(
x_{2}\right) +\lambda \right) u\left( x_{2}\right) =f\left( x_{2}\right) ,%
\text{ }L_{2j}u=0,j=1,2,
\end{equation*}%
for $x_{1}\in \left( 0,b_{1}\right) $, where $B$ is a differential operator
in $L_{p_{1}}\left( 0,b_{1};E\right) $ for $x_{2}\in \left( 0,b_{2}\right) ,$
generated by problem $\left( 2.4\right) .$ By virtue of $\left[ \text{1,
Theorem 4.5.2 }\right] $, $L_{p_{1}}\left( 0,b_{1};E\right) \in UMD$ for $%
p_{1}\in \left( 1,\infty \right) $. Moreover, in view of $\left[ 10\right] $
the operator $B$ is $R-$positive in $L_{p_{1}}\left( 0,b_{1};E\right) $.
Hence, the problem $\left( 2.5\right) $\ has a unique solution $u\in
W_{p,\alpha \left( 2\right) }^{\left[ 2\right] }\left( G_{2};E\left(
A\right) ;E\right) $ for $f\in L_{p}\left( G_{2};E\right) $, $\left\vert
\arg \lambda \right\vert \leq \varphi $ with sufficiently large $\left\vert
\lambda \right\vert $ and $\left( 2.3\right) $ holds for $n=2$. By\
continuing this process we obtain the assertion.

\bigskip \textbf{Theorem 2.2. }Let the Conditions 2.1 hold and $A_{k}\left(
x\right) A^{-\left( \frac{1}{2}-\nu \right) }\left( x\right) \in C\left( 
\bar{G};L\left( E\right) \right) $ for $0<\nu <\frac{1}{2}.$ Then, problem $%
\left( 2.1\right) $ has a unique solution $u\in W_{p\mathbf{,}\alpha }^{%
\left[ 2\right] }\left( G;E\left( A\right) ,E\right) $ for $f\in L_{p}\left(
G;E\right) $, $\left\vert \arg \lambda \right\vert \leq \varphi $ with
sufficiently large $\left\vert \lambda \right\vert $ and the coercive
uniform estimate holds

\begin{equation}
\sum\limits_{k=1}^{n}\sum\limits_{i=0}^{2}\left\vert \lambda \right\vert ^{1-%
\frac{i}{2}}\left\Vert \frac{\partial ^{\left[ i\right] }u}{\partial
x_{k}^{i}}\right\Vert _{L_{p}\left( G;E\right) }+\left\Vert Au\right\Vert
_{L_{p}\left( G;E\right) }\leq C\left\Vert f\right\Vert _{L_{p}\left(
G;E\right) }.  \tag{2.6}
\end{equation}

\textbf{Proof. }By second assumption and embedding theorem $\left[ 6\right] $
for all $h>0$ we\ have\ the following\ Ehrling-Nirenberg-Gagliardo type
estimate

\begin{equation}
\left\Vert L_{1}u\right\Vert _{L_{p}\left( G;E\right) }\leq h^{\mu }\
\left\Vert u\right\Vert _{W_{p\mathbf{,}\alpha }^{\left[ 2\right] }\left(
G;E\left( A\right) ,E\right) }+h^{-\left( 1-\mu \right) }\left\Vert
u\right\Vert _{L_{p}\left( G;E\right) }.  \tag{2.7}
\end{equation}%
Let $O$ denote the operator generated by problem $\left( 2.2\right) $ and 
\begin{equation*}
L_{1}u=\dsum\limits_{k=1}^{n}A_{k}\left( x\right) \frac{\partial ^{\left[ 1%
\right] }u}{\partial x_{k}}.
\end{equation*}%
By using the estimate $\left( 2.7\right) $ we obtain that there is a $\delta
\in \left( 0,1\right) $ such that 
\begin{equation*}
\left\Vert L_{1}\left( O+\lambda \right) ^{-1}\right\Vert _{B\left( X\right)
}<\delta .
\end{equation*}%
Hence, from perturbation theory of linear operators we obtain the assertion.

\begin{center}
\textbf{3. Abstract Cauchy problem for degenerate parabolic equation }
\end{center}

Consider the IBVP for degenerate parabolic equation with parameter:%
\begin{equation}
\frac{\partial u}{\partial t}+\dsum\limits_{k=1}^{n}a_{k}\left( x_{k}\right) 
\frac{\partial ^{\left[ 2\right] }u}{\partial x_{k}^{2}}+A\left( x\right)
u+du=f\left( x,t\right) \text{, }t\in \left( 0,T\right) \text{, }x\in G, 
\tag{3.1}
\end{equation}%
\begin{equation*}
\sum\limits_{i=0}^{m_{kj}}\alpha _{kji}u_{x_{k}}^{\left[ i\right] }\left(
G_{k0},t\right) +\beta _{kji}u_{k}^{\left[ i\right] }\left( G_{kb},t\right)
=0,\text{ }j=1,2,
\end{equation*}%
\begin{equation}
u\left( x,0\right) =0,\text{ }t\in \left( 0,T\right) \text{, }x^{\left(
k\right) }\in G_{k},  \tag{3.2}
\end{equation}%
where $u=u\left( x,t\right) $ is a solution, $\delta _{ki},$ $\beta _{ki}$
are complex numbers, $a_{k}$ are complex-valued functions on $G,$ $A\left(
x\right) $ is a linear operator in a Banach space $E$, domains $G,$ $G_{k,}$ 
$G_{k0},$ $G_{kb}$, $\sigma _{ik}$ and $x^{\left( k\right) }$\ are defined
in the section 2.

For $\mathbf{p=}\left( p_{0},p\right) ,$ $G_{T}=\left( 0,T\right) \times G,$ 
$L_{p\mathbf{,\gamma }}\left( G_{T};E\right) $ will denote the space of all $%
E$-valued weighted $\mathbf{p}$-summable functions with mixed norm.

\textbf{Theorem 3.1. }Suppose the Condition 2.1 hold for $\varphi >\frac{\pi 
}{2}$. Then, for $f\in L_{\mathbf{p}}\left( G_{T};E\right) $ and
sufficiently large $d>0$ problem $\left( 3.1\right) -\left( 3.2\right) $ has
a unique solution belonging to $W_{\mathbf{p},\alpha }^{1,\left[ 2\right]
}\left( G_{T};E\left( A\right) ,E\right) $ and the following coercive
uniform estimate holds 
\begin{equation*}
\left\Vert \frac{\partial u}{\partial t}\right\Vert _{L_{\mathbf{p}}\left(
G_{T};E\right) }+\sum\limits_{k=1}^{2}\left\Vert \frac{\partial ^{\left[ 2%
\right] }u}{\partial x_{k}^{2}}\right\Vert _{L_{\mathbf{p}}\left(
G_{T};E\right) }+\left\Vert Au\right\Vert _{L_{\mathbf{p}}\left(
G_{T};E\right) }\leq C\left\Vert f\right\Vert _{L_{\mathbf{p}}\left(
G_{T};E\right) }.
\end{equation*}%
\textbf{Proof. }The problem $\left( 3.1\right) $ can be expressed as the
following abstract Cauchy problem 
\begin{equation}
\frac{du}{dt}+\left( O+d\right) u\left( t\right) =f\left( t\right) ,\text{ }%
u\left( 0\right) =0.  \tag{3.3}
\end{equation}
By virtue of $\left[ 10\right] $,$\ O$ is $R$-positive in $X=L_{p}\left(
G;E\right) $, i.e $O$ is a generator of an analytic semigroup in $X.$ Then
by virtue of $\left[ \text{11, Theorem 4.2}\right] ,$ problem $\left(
3.3\right) $ has a unique solution $u\in W_{p_{0}}^{1}\left( 0,T;D\left(
O\right) ,X\right) $ for $f\in L_{p_{0}}\left( 0,T;X\right) $ and
sufficiently large $d>0.$ Moreover, the following uniform estimate holds 
\begin{equation*}
\left\Vert \frac{du}{dt}\right\Vert _{L_{p_{0}}\left( 0,T;X\right)
}+\left\Vert Ou\right\Vert _{L_{p_{0}}\left( 0,T;X\right) }\leq C\left\Vert
f\right\Vert _{L_{p_{0}}\left( 0,T;X\right) }.
\end{equation*}%
Since $L_{p_{0}}\left( G_{T};X\right) =L_{\mathbf{p}}\left( G_{T};E\right) ,$
by Theorem 2.1 we have 
\begin{equation*}
\left\Vert \left( O+d\right) u\right\Vert _{L_{p_{0}}\left( \left(
0,T\right) ;X\right) }=D\left( O\right) .
\end{equation*}%
Hence, the assertion follows from the above estimate.

\begin{center}
\textbf{5}. \textbf{Nonlinear degenerate abstract parabolic problem}
\end{center}

In this section, we consider IBVP for the following nonlinear degenerate
parabolic equation

\begin{equation}
\frac{\partial u}{\partial t}+\dsum\limits_{k=1}^{n}a_{k}\left( x_{k}\right) 
\frac{\partial ^{\left[ 2\right] }u}{\partial x_{k}^{2}}+B\left( \left(
t,x,u,D^{\left[ 1\right] }u\right) \right) u=F\left( t,x,u,D^{\left[ 1\right]
}u\right) ,  \tag{5.1}
\end{equation}%
\begin{equation*}
\sum\limits_{i=0}^{m_{kj}}\alpha _{kji}u_{x_{k}}^{\left[ i\right] }\left(
G_{k0},t\right) +\beta _{kji}u_{k}^{\left[ i\right] }\left( G_{kb},t\right)
=0,j=1,2,
\end{equation*}%
\begin{equation}
u\left( x,0\right) =0,\text{ }t\in \left( 0,T\right) \text{, }x\in G\text{, }%
x^{\left( k\right) }\in G_{k},  \tag{5.2}
\end{equation}

where $u=u\left( x,t\right) $ is a solution, $\alpha _{kji},$ $\beta _{kji}$
are complex numbers, $a_{k}$ are complex-valued functions on $\left[ 0,b_{k}%
\right] $; domains $G$, $G_{k}$, $G_{k0}$, $G_{kb}$ and $\sigma _{ik}$, $%
x^{\left( k\right) }$\ are defined in the section 2 and 
\begin{equation*}
D_{k}^{\left[ i\right] }u=\frac{\partial ^{\left[ i\right] }u}{\partial
x_{k}^{i}}=\left( x_{k}^{\alpha _{k}}\frac{\partial }{\partial x_{k}}\right)
^{i}u\left( x,t\right) \text{, }0\leq \alpha _{k}<1.
\end{equation*}

Let $G_{T}=\left( 0,T\right) \times G$. Moreover, we let 
\begin{equation*}
G_{0}=\prod\limits_{k=1}^{n}\left( 0,b_{0k}\right) ,\text{ }%
G=\prod\limits_{k=1}^{n}\left( 0,b_{k}\right) \text{, }b_{k}\in \left(
0,b_{0k}\right) \text{, }
\end{equation*}%
\begin{equation*}
T\in \left( 0,T_{0}\right) ,\text{ }B_{ki}=\left( W^{2,p}\left(
G_{k},E\left( A\right) ,E\right) ,L^{p}\left( G_{k};E\right) \right) _{\eta
_{ik},p},
\end{equation*}%
\begin{equation*}
\eta _{ik}=\frac{m_{ki}+\frac{1}{p\left( 1-\alpha _{k}\right) }}{2},\text{ }%
\left\vert \alpha _{kjm_{kj}}\right\vert +\left\vert \beta
_{kjm_{kj}}\right\vert >0,\text{ }B_{0}=\prod_{k=1}^{n}\prod_{i=0}^{1}B_{ki}.
\end{equation*}

Let

\begin{equation*}
\alpha =\alpha \left( x\right) =\dprod\limits_{k=1}^{n}x_{k}^{\alpha _{k}}.
\end{equation*}%
\textbf{Remark 5.0. }Under the substitutions 
\begin{equation*}
\tau _{k}=\frac{x_{k}^{1-\alpha _{k}}}{1-\alpha _{k}},\text{ }0<\alpha _{k}<1%
\text{, }k=1,2,...,n
\end{equation*}%
the spaces $L_{p}\left( G;E\right) $ and $W_{p,\alpha }^{\left[ 2\right]
}\left( G;E\left( A\right) ,E\right) $ are mapped isomorphically onto the
weighted spaces $L_{p,\tilde{\alpha}}\left( \tilde{G};E\right) $ and $W_{p,%
\tilde{\alpha}}^{2}\left( \tilde{G};E\left( A\right) ,E\right) $,
respectively, where%
\begin{equation*}
\tilde{G}=\dprod\limits_{k=1}^{n}\left( 0,\tilde{b}_{k}\right) \text{, }%
\tilde{b}_{k}=\frac{b_{k}^{1-\alpha _{k}}}{1-\alpha _{k}},\text{ }\tilde{%
\alpha}\left( \tau \right) =\alpha \left( x_{1}\left( \tau _{1}\right)
,x_{2}\left( \tau _{2}\right) ,...,x_{n}\left( \tau _{n}\right) \right) .
\end{equation*}

\textbf{Remark 5.1. }By virtue of $\left[ \text{28, \S\ 1.8.}\right] $ and
Remark 5.0, operators $u\rightarrow \frac{\partial ^{\left[ i\right] }u}{%
\partial x_{k}^{i}}\mid _{x_{k=0}}$are continuous from $W_{p,\alpha }^{\left[
2\right] }\left( G;E\left( A\right) ,E\right) $ onto $B_{ki}$ and there are
the constants $C_{1}$ and $C_{0}$ such that for $w\in W_{p,\alpha }^{\left[ 2%
\right] }\left( G;E\left( A\right) ,E\right) ,$ $W=\left\{ w_{ki}\right\} ,$ 
$w_{ki}=\frac{\partial ^{\left[ i\right] }w}{\partial x_{k}^{i}},$ $i=0,1$, $%
k=1,2,...,n$%
\begin{equation*}
\left\Vert \frac{\partial ^{\left[ i\right] }w}{\partial x_{k}^{i}}%
\right\Vert _{B_{ki},\infty }=\sup_{x\in G}\left\Vert \frac{\partial ^{\left[
i\right] }w}{\partial x_{k}^{i}}\right\Vert _{B_{ki}}\leq C_{1}\left\Vert
w\right\Vert _{W_{p,\alpha }^{\left[ 2\right] }\left( G;E\left( A\right)
,E\right) }\text{,}
\end{equation*}

\begin{equation*}
\left\Vert W\right\Vert _{0,\infty }=\sup_{x\in
G}\sum\limits_{k,i}\left\Vert w_{ki}\right\Vert _{B_{ki}}\leq
C_{0}\left\Vert w\right\Vert _{W_{p,\alpha }^{\left[ 2\right] }\left(
G;E\left( A\right) ,E\right) }.
\end{equation*}%
\textbf{Condition 5.1. }Suppose the following hold:\textbf{\ }

(1) $E$ is an UMD\ space and $0\leq \alpha _{1}$, $\alpha _{2}<1-\frac{1}{p}$%
, $p\in \left( 1,\infty \right) $;

(2) $a_{k}$ are continuous functions on $\left[ 0,b_{k}\right] ,$ $%
a_{k}\left( x_{k}\right) <0$, for all $x\in \left[ 0,b_{k}\right] ,$ $\alpha
_{km_{k1}}\neq 0,$ $\beta _{km_{k2}}\neq 0,$ $k=1,2,...,n$;

(3) there exist $\Phi _{ki}\in B_{ki}$ such that the operator $B\left(
t,x,\Phi \right) $ for $\Phi =\left\{ \Phi _{ki}\right\} \in B_{0}$ is $R$%
-positive in $E$ uniformly with respect to $x\in G_{0}$ and $t\in \left[
0,T_{0}\right] ;$ moreover, 
\begin{equation*}
B\left( t,x,\Phi \right) B^{-1}\left( t^{0},x^{0},\Phi \right) \in C\left( 
\bar{G};L\left( E\right) \right) ,\text{ }t^{0}\in \left( 0,T\right) ,\text{ 
}x^{0}\in G;
\end{equation*}

(4) $A=B\left( t^{0},x^{0},\Phi \right) $: $G_{T}\times B_{0}\rightarrow
L\left( E\left( A\right) ,E\right) $ is continuous. Moreover, for each
positive $r$ there is a positive constant $L\left( r\right) $ such that

$\left\Vert \left[ B\left( t,x,U\right) -B\left( t,x,\bar{U}\right) \right]
\upsilon \right\Vert _{E}\leq L\left( r\right) \left\Vert U-\bar{U}%
\right\Vert _{B_{0}}\left\Vert A\upsilon \right\Vert _{E}$

for $t\in \left( 0,T\right) ,$ $x\in G$, $U,\bar{U}\in B_{0},\bar{U}=\left\{ 
\bar{u}_{ki}\right\} ,\bar{u}_{ki}\in B_{ki},$ $\left\Vert U\right\Vert
_{B_{0}},\left\Vert \bar{U}\right\Vert _{B_{0}}\leq r,\upsilon \in D\left(
A\right) ;$

(5) the function $F$: $G_{T}\times B_{0}\rightarrow E$ such that $F\left(
.,U\right) $ is measurable for each $U\in B_{0}$ and $F\left( t,x,.\right) $
is continuous for a.a. $t\in \left( 0,T\right) $, $x\in G.$ Moreover, $%
\left\Vert F\left( t,x,U\right) -F\left( t,x,\bar{U}\right) \right\Vert
_{E}\leq \Psi _{r}\left( x\right) \left\Vert U-\bar{U}\right\Vert _{B_{0}}$
for a.a. $t\in \left( 0,T\right) $, $x\in G$, $U,\bar{U}\in B_{0}$ and $%
\left\Vert U\right\Vert _{B_{0}},\left\Vert \bar{U}\right\Vert _{B_{0}}\leq
r $; $f\left( .\right) =F\left( .,0\right) \in L_{p}\left( G_{T};E\right) $.

The main result of this section is the following:

\textbf{Theorem 5.1. }Let Condition 5.1 hold. Then there is a $T\in \left(
0,T_{0}\right) $ and a $b_{k}\in \left( 0,b_{0k}\right) $ such that\ problem 
$\left( 5.1\right) -\left( 5.2\right) $ has a unique solution belonging to $%
W_{p,\alpha }^{1,\left[ 2\right] }\left( G_{T};E\left( A\right) ,E\right) .$

\textbf{Proof. }Consider the following linear problem%
\begin{equation*}
\frac{\partial w}{\partial t}+\dsum\limits_{k=1}^{n}a_{k}\left( x_{k}\right) 
\frac{\partial ^{\left[ 2\right] }w}{\partial x_{k}^{2}}+du=f\left(
x,t\right) ,\text{ }x\in G,\text{ }t\in \left( 0,T\right) ,
\end{equation*}%
\begin{equation}
\sum\limits_{i=0}^{m_{kj}}\alpha _{kji}w_{x_{k}}^{\left[ i\right] }\left(
G_{k0},t\right) +\sum\limits_{i=0}^{m_{k2}}\beta _{kji}w_{k}^{\left[ i\right]
}\left( G_{kb},t\right) =0,j=1,2,  \tag{5.3}
\end{equation}%
\begin{equation*}
w\left( x,0\right) =0,t\in \left( 0,T\right) \text{, }x\in G\text{, }%
x^{\left( k\right) }\in G_{k},\text{ }d>0.
\end{equation*}

By Theorem 3.1and in view of Proposition 4.1 there is a unique solution $%
w\in W_{p,\alpha }^{1,\left[ 2\right] }\left( G_{T};E\left( A\right)
,E\right) $ of the problem $(5.3)$ for $f\in L_{p}\left( G_{T};E\right) $
and sufficiently large $d>0$ and it satisfies the following coercive
estimate 
\begin{equation*}
\left\Vert w\right\Vert _{W_{p,\alpha }^{1,\left[ 2\right] }\left(
G_{T};E\left( A\right) ,E\right) }\leq C_{0}\left\Vert f\right\Vert
_{L_{p}\left( G_{T};E\right) },
\end{equation*}%
\ uniformly with respect to $b\in \left( 0\right. ,\left. b_{0}\right] $,
i.e., the constant $C_{0}$\ does not depends on $f\in L_{p}\left(
G_{T};E\right) $ and $b\in \left( 0\right. \left. b_{0}\right] $ where 
\begin{equation*}
A\left( x\right) =B\left( x,0\right) ,\text{ }f\left( x\right) =F\left(
x,0\right) ,\text{ }x\in \left( 0,b\right) .
\end{equation*}%
We want to solve the problem $\left( 5.1\right) -\left( 5.2\right) $ locally
by means of maximal regularity of the linear problem $(5.3)$ via the
contraction mapping theorem. For this purpose, let $w$ be a solution of the
linear BVP $(5.3).$ Consider a ball 
\begin{equation*}
B_{r}=\left\{ \upsilon \in Y,\text{ }\upsilon -w\in Y_{1},\text{ }\left\Vert
\upsilon -w\right\Vert _{Y}\leq r\right\} .
\end{equation*}

For given $\upsilon \in B_{r},$ consider the following linearized problem%
\begin{equation}
\frac{\partial u}{\partial t}+\dsum\limits_{k=1}^{n}a_{k}\left( x_{k}\right) 
\frac{\partial ^{\left[ 2\right] }u}{\partial x_{k}^{2}}+A\left( x\right)
=F\left( x,V\right) +\left[ B\left( x,0\right) -B\left( x,V\right) \right]
\upsilon ,  \notag
\end{equation}%
\begin{equation}
\sum\limits_{i=0}^{m_{kj}}\alpha _{kji}w_{x_{k}}^{\left[ i\right] }\left(
G_{k0},t\right) +\sum\limits_{i=0}^{m_{k2}}\beta _{kji}w_{x_{k}}^{\left[ i%
\right] }\left( G_{kb},t\right) =0,  \tag{5.4}
\end{equation}%
\begin{equation*}
w\left( x,0\right) =0,\text{ }t\in \left( 0,T\right) \text{, }x\in G\text{, }%
x^{\left( k\right) }\in G_{k}.
\end{equation*}%
where $V=\left\{ \upsilon _{ki}\right\} ,$ $\upsilon _{ki}\in B_{ki}.$
Define a map $Q$ on $B_{r}$ by $Q\upsilon =u,$ where $u$ is solution of $%
\left( 5.4\right) .$ We want to show that $Q\left( B_{r}\right) \subset
B_{r} $ and that $Q$ is a contraction operator provided $T$ and $b_{k}$ are
sufficiently small, and $r$ is chosen properly. In view of separability
properties of the problem $\left( 5.3\right) $ we have%
\begin{equation*}
\left\Vert Q\upsilon -w\right\Vert _{Y}=\left\Vert u-w\right\Vert _{Y}\leq
C_{0}\left\{ \left\Vert F\left( x,V\right) -F\left( x,0\right) \right\Vert
_{X}+\right.
\end{equation*}

\begin{equation*}
\left. \left\Vert \left[ B\left( 0,W\right) -B\left( x,V\right) \right]
\upsilon \right\Vert _{X}\right\} .
\end{equation*}

By assumption (4) we have 
\begin{equation*}
\left\Vert \left[ B\left( 0,W\right) \upsilon -B\left( x,V\right) \right]
\upsilon \right\Vert _{X}\leq \sup\limits_{x\in \left[ 0,b\right] }\left\{
\left\Vert \left[ B\left( 0,W\right) -B\left( x,W\right) \right] \upsilon
\right\Vert _{L\left( E_{0},E\right) }\right.
\end{equation*}%
\begin{equation*}
+\left. \left\Vert B\left( x,W\right) -B\left( x,V\right) \right\Vert
_{L\left( E_{0},E\right) }\left\Vert \upsilon \right\Vert _{Y}\right\} \leq
\end{equation*}%
\begin{equation*}
\left[ \delta \left( b\right) +L\left( R\right) \left\Vert W-V\right\Vert
_{\infty ,E_{0}}\right] \left[ \left\Vert \upsilon -w\right\Vert
_{Y}+\left\Vert w\right\Vert _{Y}\right] \leq
\end{equation*}%
\begin{equation*}
\left\{ \delta \left( b\right) +L\left( R\right) \left[ C_{1}\left\Vert
\upsilon -w\right\Vert _{Y}+\left\Vert \upsilon -w\right\Vert _{Y}\right]
\right.
\end{equation*}%
\begin{equation*}
\left. \left[ \left\Vert \upsilon -w\right\Vert _{Y}+\left\Vert w\right\Vert
_{Y}\right] \right\} \leq \delta \left( b\right) +L\left( R\right) \left[
C_{1}r+r\right] \left[ r+\left\Vert w\right\Vert _{Y}\right] ,
\end{equation*}

where 
\begin{equation*}
\delta \left( b\right) =\sup\limits_{x\in \left[ 0,b\right] }\left\Vert %
\left[ B\left( 0,W\right) -B\left( x,W\right) \right] \right\Vert _{B\left(
E_{0},E\right) }.
\end{equation*}

By assumption (5) we get 
\begin{equation*}
\left\Vert F\left( x,V\right) -F\left( x,0,\right) \right\Vert _{E}\leq
\delta \left( b\right) +
\end{equation*}%
\begin{equation*}
\left\Vert F\left( x,V\right) -F\left( x,W\right) \right\Vert
_{E}+\left\Vert F\left( x,W\right) -F\left( x,0\right) \right\Vert _{E}\leq
\end{equation*}

\begin{equation*}
\delta \left( b\right) +\mu _{R}\left[ \left\Vert \upsilon -w\right\Vert
_{Y}+\left\Vert w\right\Vert _{Y}\right]
\end{equation*}%
\begin{equation*}
\mu _{R}C_{1}\left[ \left\Vert \upsilon -w\right\Vert _{Y}+\left\Vert
w\right\Vert _{Y}\right] \leq \mu _{R}\left[ C_{1}r+\left\Vert w\right\Vert
_{Y}\right] ,
\end{equation*}%
where $R=C_{1}r+\left\Vert w\right\Vert _{Y}$ is a fixed number. In view of
above estimates, by suitable choice of $\mu _{R},$ $L_{R}$ and for
sufficiently small $T\in \left( 0,T_{0}\right) $ and $b_{k}\in \left(
0\right. ,\left. b_{0k}\right] $ we have 
\begin{equation*}
\left\Vert Q\upsilon -w\right\Vert _{Y}\leq r,
\end{equation*}%
i.e. 
\begin{equation*}
Q\left( B_{r}\right) \subset B_{r}.
\end{equation*}%
Moreover, in a similar way we obtain 
\begin{equation*}
\left\Vert Q\upsilon -Q\bar{\upsilon}\right\Vert _{Y}\leq C_{0}\left\{ \mu
_{R}C_{1}+M_{a}+L\left( R\right) \left[ \left\Vert \upsilon -w\right\Vert
_{Y}+C_{1}r\right] \right. +
\end{equation*}%
\begin{equation*}
\left. L\left( R\right) C_{1}\left[ r+\left\Vert w\right\Vert _{Y}\right]
\left\Vert \upsilon -\bar{\upsilon}\right\Vert _{Y}\right\} +\delta \left(
b\right) .
\end{equation*}

By suitable choice of $\mu _{R},$ $L_{R}$\ and for sufficiently small $T\in
\left( 0,T_{0}\right) $ and $b_{k}\in \left( 0,b_{0k}\right) $ we obtain $%
\left\Vert Q\upsilon -Q\bar{\upsilon}\right\Vert _{Y}<\eta \left\Vert
\upsilon -\bar{\upsilon}\right\Vert _{Y},$ $\eta <1,$ i.e. $Q$ is a
contraction operator. Eventually, the contraction mapping principle implies
a unique fixed point of $Q$ in $B_{r}$ which is the unique strong solution $%
u\in W_{p,\alpha }^{1,\left[ 2\right] }\left( G_{T};E\left( A\right)
,E\right) .$\textbf{\ }

\bigskip\ \textbf{References}

\begin{quote}
\ \ \ \ \ \ \ \ \ \ \ \ \ \ \ \ \ \ \ \ \ \ \ \ \ \ \ \ \ \ \ \ \ \ \ \ \ \
\ \ \ \ \ \ \ \ \ \ \ \ \ \ \ \ \ \ \ \ \ \ \ \ \ \ \ \ \ \ \ \ \ \ \ \ \ \
\ \ \ \ \ \ \ \ \ \ \ \ \ \ \ \ \ 
\end{quote}

\begin{enumerate}
\item H. Amann, Linear and quasilinear parabolic problems,1, 2, Birkhauser,
Basel 1995, MR134538

\item S. Yakubov and Ya. Yakubov, Differential-operator Equations. Ordinary
and Partial \ Differential Equations, Chapman and Hall /CRC, Boca Raton,
2000, Zbl 0936.35002

\item A. Favini, A. Yagi, Degenerate Differential Equations in Banach
Spaces, Taylor \& Francis, Dekker, New-York, 1999, MR1654663, Zbl 0913.34001

\item R. Agarwal R, D. O' Regan, V. B. Shakhmurov, Separable anisotropic
differential operators in weighted abstract spaces and applications, J.
Math. Anal. Appl. (338(2008), 970-983, MR2386473,
DOI:10.1016/j.jmaa.2007.05.078

\item C. Ashyralyev, C. Claudio and S. Piskarev, On well-posedness of
difference schemes for abstract elliptic problems in $L_{p}$ spaces, Numer.
Func. Anal.Optim., (29)1, 2 (2008), 43-65, MR2387837,
DOI:10.1080/01630560701872698

\item V. B. Shakhmurov, Embedding theorems and maximal regular differential
operator equations in Banach-valued function spaces, J. Inequal. Appl., (4)
2005, 605-620, MR2210722; DOI: 10.1155/JIA.2005.329

\item V. B. Shakhmurov, Linear and nonlinear abstract equations with
parameters, Nonlinear Anal., 73(2010), 2383-2397, MR2674077;
DOI:10.1016/j.na.2010.06.004

\item V. B. Shakhmurov, Regular degenerate separable differential operators
and applications, Potential Anal., 35(3) (2011), 201-212, MTR2832575; DOI:
10.1007/s11118-010-9206-9

\item V. B. Shakhmurov, Coercive boundary value problems for regular
degenerate differential-operator equations, J. Math. Anal. Appl., 292 (2),
(2004), 605-620, MR2048274; DOI:10.1016/j.jmaa.2003.12.032

\item V. B., Shakhmurov, A. Shahmurova, Nonlinear abstract boundary value
problems atmospheric dispersion of pollutants, Nonlinear Anal. Real World
Appl., 11 (2) (2010), 932-951; MR2571266, DOI:10.1016/j.nonrwa.2009.01.037

\item L. Weis, Operator-valued Fourier multiplier theorems and maximal $%
L_{p} $ regularity, Math. Ann. 319(2001), 735-758, MR1825406; DOI:
10.1007/PL00004457

\item W. E. Fitzgibbon, M. Langlais, J. J. Morgan, A degenerate
reaction-diffusion system modeling atmospheric dispersion of pollutants, J.
Math. Anal. Appll. \ 307(2005), 415-432, MR2142434; DOI:
10.1016/j.jmaa.2005.02.060

\item D. L. Burkholder, A geometrical characterization of Banach spaces in
which martingale difference sequences are unconditional, Ann. Probab. 9(6)
(1981), 997--1011, MR 632972, DOI:10.1214/aop/1176994270
\end{enumerate}

\end{document}